\documentclass[12pt,reqno]{amsart}
\usepackage{enumerate, latexsym, amsmath, amsfonts, amssymb, amsthm, color}
\usepackage{psfrag,graphicx, colordvi}
\def\pmod #1{\ ({\rm{mod}}\ #1)}
\def\Z{\mathbb Z}

\def\Q{\mathbb Q}

\def\l{\left}
\def\r{\right}
\def\bg{\bigg}
\def\({\bg(}
\def\){\bg)}
\def\t{\text}
\def\f{\frac}

\def\ls{\leqslant}

\def\al{\alpha}

\def\ve{\varepsilon}

\def\eq{\equiv}

\def\da{\delta}

\def\Proof{\noindent{\it Proof}}

\theoremstyle{plain}
\newtheorem{theorem}{Theorem}

\newtheorem{lemma}{Lemma}

\theoremstyle{definition}

\theoremstyle{remark}
\newtheorem{remark}{Remark}

\makeatletter
\@namedef{subjclassname@2020}{%
  \textup{2020} Mathematics Subject Classification}
\makeatother
 \vspace{4mm}

\makeatletter \@addtoreset{equation}{section}
\@addtoreset{theo}{section} \makeatother

\def\qed{\hfill \rule{4pt}{7pt}}

\begin{document}
\baselineskip=17pt
\medskip

\hbox{C. R. Math. Acad. Sci. Paris 360 (2022), 1065--1069.}
\medskip

\title
[A new theorem on quadratic residues modulo primes]
{A new theorem on quadratic residues modulo primes}

\author
[Qing-Hu Hou, Hao Pan and Zhi-Wei Sun]
{Qing-Hu Hou, Hao Pan and Zhi-Wei Sun}

 \address {(Qing-Hu Hou) School of Mathematics, Tianjin University,
  Tianjin 300350, People's Republic of China}

\email{{\tt qh\_hou@tju.edu.cn}}

\address {(Hao Pan) School of Applied Mathematics, Nanjing University of Finance and Economics,
 Nanjing 210046, People's Republic of China}

\email {{\tt haopan79@zoho.com}}

\address {(Zhi-Wei Sun, corresponding author) Department of Mathematics, Nanjing
University, Nanjing 210093, People's Republic of China}

\email{{\tt zwsun@nju.edu.cn}
\newline\indent
{\it Homepage}: {\tt http://maths.nju.edu.cn/\lower0.5ex\hbox{\~{}}zwsun}}


\thanks{The first, second and third authors are supported by the
  National Natural Science Foundation of China (grants 11771330, 12071208 and 11971222,
  respectively)}

\begin{abstract} Let $p>3$ be a prime, and let $(\frac{\cdot}p)$ be the Legendre symbol.
Let $b\in\mathbb Z$ and $\varepsilon\in\{\pm1\}$.
We mainly prove that
$$\left|\left\{N_p(a,b):\ 1<a<p\ \text{and}\ \left(\frac ap\right)=\varepsilon\right\}\right|=\frac{3-(\frac{-1}p)}2,$$
where $N_p(a,b)$ is the number of positive integers $x<p/2$ with $\{x^2+b\}_p>\{ax^2+b\}_p$, and
$\{m\}_p$ with $m\in\mathbb Z$ is the least nonnegative residue of $m$ modulo $p$.
\end{abstract}

\maketitle

\section{Introduction}

The theory of quadratic residues modulo primes plays an important role in fundamental number theory.

Let $p$ be an odd prime and let $a\in\Z$ with $p\nmid a$.
By Gauss' Lemma (cf. \cite[p.\,52]{IR}),
$$\l(\f ap\r)=(-1)^{|\{1 \ls k\ls\f{p-1}2:\ \{ka\}_p>\f p2\}|},$$
where $(\f{\cdot}p)$ denotes the Legendre symbol, and we write $\{x\}_p$ for the least nonnegative residue of an integer $x$ modulo $p$.

Let $n$ be any positive odd integer, and let $a\in\Z$ with $\gcd(a(1-a),n)=1$. In 2020, Z.-W. Sun
\cite{SunIJNT} proved the following new result:
$$(-1)^{|\{1\ls k\ls\f{n-1}2:\ \{ka\}_n>k\}|}=\l(\f{2a(1-a)}n\r),$$
where $(\f{\cdot}n)$ is the Jacobi symbol.

Let $p$ be an odd prime and let $a,b\in\Z$ with $a(1-a)\not\eq0\pmod p$. By \cite[Lemma 2.7]{S19}, we have
$$|\{x\in\{0,\ldots,p-1\}:\ \{ax+b\}_p>x\}|=\f{p-1}2.$$
In 2019 Z.-W. Sun \cite{S19} employed Galois theory to prove that
$$(-1)^{|\{1\ls i<j\ls \f{p-1}2:\ \{i^2\}_p>\{j^2\}_p\}|}=\begin{cases}1&\t{if}\ p\eq3\pmod8,
\\(-1)^{(h(-p)+1)/2}&\t{if}\ p\eq7\pmod8,\end{cases}$$
where $h(-p)$ is the class number of the imaginary quadratic field $\Q(\sqrt{-p})$.

Motivated by the above work, for an odd prime $p$ and integers $a$ and $b$, we introduce the notation
$$N_p(a,b):=\l|\l\{1\ls x\ls\f{p-1}2:\ \{x^2+b\}_p>\{ax^2+b\}_p\r\}\r|.$$

{\it Example} 1.1. We have $N_7(4,0)=2$ since
$$\{1^2\}_7<\{4\times1^2\}_7,\ \{2^2\}_7>\{4\times2^2\}_7\ \t{and}\ \{3^2\}_7>\{4\times3^2\}_7.$$

Let $p$ be a prime with $p\eq1\pmod4$.
 Then $q^2 \eq-1\pmod p$ for some integer $q$, hence for $a,x\in\Z$ we have
$\{(qx)^2\}_p > \{a(qx)^2\}_p$ if and only if $\{x^2\}_p < \{ax^2\}_p.$ Thus, for each $a = 2,\ldots,p-1$ there are exactly $(p-1)/4$ positive integers $x < p/2$ such that $\{x^2\}_p>\{ax^2\}_p$. Therefore
$N_p(a,0) = (p-1)/4$ for all $a = 2,\ldots,p-1$.

In this paper we establish the following novel theorem which was conjectured by the first and third authors \cite{HS} in 2018.

\begin{theorem} \label{Main} Let $p > 3$ be a prime, and let $b$ be any integer. Set
$$S = \l\{N_p(a,b):\ 1 < a < p\ \t{and}\ \l(\f ap\r)=1\r\}$$
and
 $$T = \l\{N_p(a,b):\ 1 < a < p\ \t{and}\ \l(\f ap\r)=-1\r\}.$$
Then $|S|=|T|=1$ if $p\eq1\pmod4$, and $|S|=|T|=2$ if $p\eq3\pmod 4$.
 Moreover, the set $S$ does not depend on the value of $b$.
 \end{theorem}

{\it Example 1.2}. Let's adopt the notation in Theorem \ref{Main}. For $p=5$, we have $S=\{1\}$ for any $b \in \mathbb{Z}$, and the set $T$ depends on $b$ as illustrated by the following table:
\begin{center}
\begin{tabular}{|c|c|c|c|c|c|} \hline
$b$&\  0\ &\ 1\ &\ 2\ &\
3\ &\ 4\ \\
\hline $T$&\{1\}&\{0\}&\{1\}&\{2\}&\{1\}\\
\hline
\end{tabular}.
\end{center}
For $p=7$, we have $S=\{1,2\}$ for any $b \in \mathbb{Z}$,  and the set $T$ depends on $b$ as illustrated by the following table:
\begin{center}
\begin{tabular}{|c|c|c|c|c|c|c|c|}\hline
$b$ & 0 & 1 & 2 & 3 & 4 & 5 & 6  \\ \hline
\rule{0pt}{15pt}
$T$ & \{0,1\} & \{1,2\} & \{2,3\} & \{1,2\} & \{2,3\} & \{1,2\} & \{0,1\}\\
\hline
\end{tabular}.\end{center}

\section{Proof of Theorem 1.1}
\label{sec:2}

\begin{lemma}\label{Lem2.1} For any prime $p\eq3\pmod4$, we have
 $$\sum_{z=1}^{p-1}z\l(\f zp\r)=-ph(-p),$$
 where $h(-p)$ is the class number of the imaginary quadratic field $\Q(\sqrt{-p})$.
\end{lemma}
\begin{remark} This is a known result of Dirichlet (cf. \cite[Corollary 5.3.13]{Cohen93}).
\end{remark}

\begin{lemma}\label{Lem2.2} For any prime $p\eq3\pmod4$ with $p>3$, there are $x,y,z\in\{1,\ldots,p-1\}$
such that $$\l(\f xp\r)=\l(\f {x+1}p\r)=1, \ -\l(\f yp\r)=\l(\f{y+1}p\r)=1,
\ \t{and}\ \l(\f zp\r)=-\l(\f{z+1}p\r)=1.$$
\end{lemma}
\Proof. By a known result (see, e.g., \cite[pp.\,64--65]{D}), we have
$$\l|\l\{x\in\{1,\ldots,p-2\}:\ \l(\f xp\r)=\l(\f{x+1}p\r)=1\r\}\r|=\f{p-3}4>0.$$
Hence
\begin{align*}&\l|\l\{y\in\{1,\ldots,p-2\}:\ -\l(\f yp\r)=\l(\f{y+1}p\r)=1\r\}\r|
\\=&\l|\l\{y\in\{1,\ldots,p-2\}:\ \l(\f{y+1}p\r)=1\r\}\r|-\f{p-3}4
\\=&\f{p-1}2-1-\f{p-3}4=\f{p-3}4>0
\end{align*}
and
\begin{align*}&\l|\l\{z\in\{1,\ldots,p-2\}:\ \l(\f zp\r)=-\l(\f{z+1}p\r)=1\r\}\r|
\\=&\l|\l\{z\in\{1,\ldots,p-2\}:\ \l(\f{z}p\r)=1\r\}\r|-\f{p-3}4
\\=&\f{p-1}2-\f{p-3}4=\f{p+1}4>0.
\end{align*}
Now the desired result immediately follows. \qed

\medskip
\noindent{\bf Proof of Theorem 1.1}.  Let $a\in\{2,\ldots,p-1\}$. For any $x\in\Z$, it is easy to see that
 $$\l\{\f{ax^2+b}p\r\}+\l\{\f{(1-a)x^2}p\r\}-\l\{\f{x^2+b}p\r\}
 =\begin{cases}0&\t{if}\ \{x^2+b\}_p>\{ax^2+b\}_p,
 \\1&\t{if}\ \{x^2+b\}_p<\{ax^2+b\}_p,\end{cases}$$
 where $\{\al\}$ denotes the fractional part of a real number $\al$.
 Thus
 \begin{align*}N_p(a,b)=&\sum_{x=1}^{(p-1)/2}\l(1+\l\{\f{x^2+b}p\r\}-\l\{\f{ax^2+b}p\r\}
 -\l\{\f{(1-a)x^2}p\r\}\r)
 \\=&\f{p-1}2+\sum_{x=1}^{(p-1)/2}\l\{\f{x^2+b}p\r\}-\sum_{x=1}^{(p-1)/2}\l\{\f{ax^2+b}p\r\}
 -\sum_{x=1}^{(p-1)/2}\l\{\f{(1-a)x^2}p\r\}
 \\=&\f{p-1}2+\sum_{x=1\atop (\f xp)=1}^{p-1}\l\{\f{x+b}p\r\}
 -\sum_{y=1\atop (\f yp)=(\f ap)}^{p-1}\l\{\f{y+b}p\r\}
 -\sum_{z=1\atop (\f zp)=(\f {1-a}p)}^{p-1}\f{z}p.
 \end{align*}
 Suppose that $(\f ap)=\ve$ with $\ve\in\{\pm1\}$. Then
 $$N_p(a,b)=\f{p-1}2+\sum_{x=1\atop (\f xp)=1}^{p-1}\l\{\f{x+b}p\r\}
 -\sum_{y=1\atop (\f yp)=\ve}^{p-1}\l\{\f{y+b}p\r\}
 -\sum_{z=1\atop (\f zp)=\da\ve}^{p-1}\f{z}p,$$
 where $\da=(\f{a(1-a)}p)$.

 If $\ve=1$, then
 $$N_p(a,b)=\f{p-1}2-\f1p\sum_{z=1\atop (\f zp)=\da}^{p-1}z$$
 does not depend on $b$.

 If $p\eq1\pmod 4$, then $(\f{-1}p)=1$ and hence
 $$\sum_{z=1\atop (\f zp)=1}^{p-1}z =\sum_{z=1\atop(\f {p-z}p)=1}^{p-1}(p-z)
 =p\f{p-1}2-\sum_{z=1\atop(\f zp)=1}^{p-1}z,$$
 thus
 $$\sum_{z=1\atop (\f zp)=1}^{p-1}z =p\f{p-1}4$$
 and
 $$\sum_{z=1\atop (\f zp)=-1}^{p-1}z=\sum_{z=1}^{p-1}z-p\f{p-1}4=p\f{p-1}4.$$
 So, if $p\eq1\pmod4$, then  $|S|=|T|=1$, and moreover
 $$S=\l\{\f{p-1}2-\f{p-1}4\r\}=\l\{\f{p-1}4\r\}.$$

 Now assume that $p\eq3\pmod4$. We want to show that $|S|=|T|=2$. By Lemma \ref{Lem2.1},
 $$\sum_{z=1}^{p-1}z\l(\f zp\r)=-ph(-p)\not=0.$$
 Thus
 $$\sum_{z=1\atop (\f zp)=1}^{p-1}z=\sum_{z=1}^{p-1}z\f{1+(\f zp)}2=p\f{p-1}4-\f p2h(-p)$$
 and hence
 $$\sum_{z=1\atop (\f zp)=-1}^{p-1}z=\sum_{z=1}^{p-1}z-\sum_{z=1\atop (\f zp)=1}^{p-1}z
 =p\f{p-1}4+\f p2h(-p).$$
 By Lemma \ref{Lem2.2}, for some $a\in\{2,\ldots,p-2\}$ we have $(\f{a-1}p)=(\f ap)=1$
 and hence $(\f{a(1-a)}p)=-1$. For $a'=p+1-a$, we have
 $$\l(\f{a'}p\r)=-1\ \t{and}\ \l(\f{a'(1-a')}p\r)=\l(\f{(1-a)a}p\r)=-1.$$
 By Lemma \ref{Lem2.2}, for some $a_*,b_*\in\{2,\ldots,p-2\}$ we have
 $$-\l(\f{a_*-1}p\r)=\l(\f {a_*}p\r)=1\ \ \t{and}\ \ \l(\f{b_*-1}p\r)=-\l(\f {b_*}p\r)=1.$$
Note that
$$\l(\f{a_*(1-a_*)}p\r)=1=\l(\f{b_*(1-b_*)}p\r).$$
 Now we clearly have $|S|=|T|=2$. Moreover,
 $$S=\l\{\f{p-1}2-\l(\f{p-1}4\pm\f{h(-p)}2\r)\r\}=\l\{\f{p-1\pm2h(-p)}4\r\}.$$

The proof of Theorem 1.1 is now complete. \qed

\end{document}